\documentclass[12pt, twoside, final]{article}

\usepackage{latexsym}

\newtheorem{Th}{Theorem}
\newtheorem{Prop}{Proposition}

\newtheorem{Hup}[Th]{Conjecture}
\newtheorem{Def}[Prop]{Definition}

\begin{document}

\title{Several Remarks on Pascal Automorphism and Infinite Ergodic Theory}

\maketitle


\begin{center}
	{\large A.~M.~Vershik}\footnote{The research is supported by the Russian Science Foundation grant
		14-11-00581.}\\[3pt]
	{
		St.~Petersburg Department of Steklov Institute of Mathematics\\
		Mathematical Department of. St.~Petersburg State University\\
		Moscow Institute for Information  Transmission Problems
	}
\end{center}

{
	\hfill{
	\small \em	To Professor Arshag Hajian in connection with his jubilee}
}

\begin{abstract}
 We interpret the Pascal-adic transformation as a generalized induced automorphism (over odometer) and formulate the   $\sigma$-finite analog of odometer which is also known as "Hajian-Kakutani transformation" (former  "Ohio state example").  We shortly suggest a sketch of the theory of random walks on the groups on the base of   $\sigma$-finite ergodic theory.
\end{abstract}

\makeatletter
 \renewcommand{\@evenhead}{\scriptsize {{\thepage}
 \centerline {A.M. Vershik}}}
 \renewcommand{\@oddhead}{\scriptsize {\centerline
{Several Remarks on Pascal Automorphism and Infinite Ergodic Theory}} \hfil{\thepage}}

\bigskip

Infinite Ergodic theory, i.e. the abstract theory of transformations with an infinite ($\sigma$-finite invariant measure)
began with the well-known work of E.~Hopf and then became one of the  branches of Ergodic theory.
Many of the concepts and facts of this theory failed to transfer on the case of actions with infinite measure,
while not all of them were automatic.
However, the authentic specificity of the area was revealed after the discovery by S.~Kakutani and his disciple A.~Hajian so-called $eww$-sets.
In my opinion, the theory of these sets are still not  took the worthy place in the dynamics. In this note,
I do not touch this issue, and write only about the problems more close to me, which also are  relevant to these sets.
I pay tribute to Prof. Hajian for his long-standing efforts in this area and for his  faithfulness to the selected topic.

\section{Short history of Pascal automorphism} Pascal automorphism was used (without this name) in the paper by S.~Kakutani \cite{Ka} and coauthors in connection with a beautiful combinatorial problem about equidistribution of sequences of partitions, and also in connection with $\sigma$-finite ergodic theory.
 The main observation was in the paper by Kakutani \cite{Ka} where
 he considered that combinatorial problem:
 {\it the formulas (8),(9) on the page 266 in \cite{Ka} is just
 the formula of Pascal automorphism as measure preserving transformation of the unit interval with Lebesgue measure. It one-to-one coincides with my formula of the next paragraph which
 I suggested in 1981 \cite{V82}.}

 In the paper by A.~Hajan, Y.~Ito and S.~Kakutani \cite{HIK} it was proved that (Pascal) automorphism is ergodic (with respect to any Bernoulli measure on the interval). It is the key moment
 in the proof of equidistribution in the Kakutani Problem. But in \cite{HIK} Pascal automorphism was used for needs of $\sigma$-finite
 ergodic theory. Namely the integral model over Pascal automorphism
 was later called "Ohio-State example" (OSE) and now "Hajian-Kakutani transformation" after one of the first presentations in
 \cite{AK}, and later was developed by A.~Hajian and his collaborators.
  This is a natural example of the big series of transformations with
  $\sigma$-finite invariant measure.

  I did not know about papers \cite{Ka,HIK} before my visit to NE-University in 2011, when during my talk about Pascal automorphisms on seminar of Prof. A.~Hajian the participants  pointed out that Pascal automorphism they knew but from a different point of view.

 I defined the Pascal automorphism in the paper \cite{V82} as a nice and simplest non-trivial case of what I had called "adic transformation". \footnote{"adic" means "p-adic"  without "p",
 like Arnold's notion of "versal" deformation means "universal" without "uni-".}
 Adic transformations was defined in the end of 70-th (see {\cite{V81}) in order to develop the strong approximation in ergodic theory. Adic realization of the transformation
 is nothing more than sequence of the coherent Rokhlin towers
 considered in the space of paths of the graph, and what is important the height of towers can depend on the point and can change with $n$.
 \footnote{The Dye's theorem about isomorphism of orbit partitions of any ergodic transformations is an easy corollary
 of adic realization (see \cite{V68}). The analog of Dye's theorem for ergodic transformation with $\sigma$-finite invariant measure  (\cite{Krie}) is corollary of the integral realization of automorphism (see next paragraph).}

 Now adic transformations became very popular --- we can speak about adic-type of dynamics as a theory of the special type of dynamical systems like symbolic dynamics. This is a sort of constructions of important examples and counter-examples in the theory. Shortly speaking adic dynamics is dynamics of paths
 of the graded graphs (Bratteli diagrams) rather action of the group.
 \footnote{The popular term "Bratteli-Vershik diagram" is not precise: the "adic transformation" was defined by myself as a model of approximation in the uniform topology of  automorphisms,
 and were a self-consistent collections of Rokhlin towers for given automorphisms. The interpretation with Bratteli diagram was very
 convenient if we define the lexicographic order on the paths.
 The corresponding automorphisms and their - names for the given graphs (Pascal, Young etc.) I had used in the beginning of 80-th.}

  From topological point of view this is a dynamics on Cantor set or more exactly on the space of paths of a graded graph. From symbolic point
 of view this is a "transversal dynamics" on the Markov or more general but in general not stationary compact. For the stationary case this definition had been considered by S.Ito \cite{Ito}. Measure-theoretical and Borel approach to the adic dynamics is useful for many asymptotical and probabilistic
 problems about combinatorics and topological structure of space of paths or Markov trajectories. As the examples we can mention limit shape theorems, approximation, entropy and so on, see my survey \cite{VJa}. In this notes I will mention several questions which related to
Pascal automorphism as well as to the $\sigma$-finite Ergodic theory.
We do not touch here more deeper properties of automorphisms
with infinite invariant measures, like {\em ww}-sets and others.

  In the next section we made some definitions and define a new notion of generalized induced automorphisms which shows that Pascal
  automorphism is induced in the sense of new definition of odometer.

  In the last section we gave a sketch of the $\sigma$-finite version
  of the theory of random walks on  the groups.

 \section{The odometer and Pascal automorphism}

 Recall the definitions of odometer and Pascal automorphisms.
 Consider $X=\textbf{Z}_2$, the  compact additive dyadic group of
dyadic integers with Haar measure $\mu$ and let  $$T:Tx=x+1$$
be the addition of unity, or odometer-ergodic transformation
with dyadic spectra. The $p$-adic odometer can be define in the same way.
Moreover we can consider an arbitrary product-space of the finite
sets
 $$\prod_k Z_{p_k},\,  p_k \in ({\textbf N}\setminus 1),\,  k=1,2 \dots,$$
 equipped with group structure as profinite group, and then define
 $\{p_k\}$-odometer. We will consider here only 2- or $p$-odometer.

    Denote as $P$ {\it the Pascal  automorphism of the space $(Z_2, \mu)$} an automorphism which is  defined by the  following
   formula in terms of dyadic decomposition:
  $$x\mapsto Px;\quad
  P(0^m1^k\textbf{10}**),=1^k0^m\textbf{01}** \quad m,k=0,1
  \dots.$$
   For example: $$P(10**)=01**, P(0110**)=1001,
   P(00110**)=10001,
   P(1110**)=1101**,$$ etc.
   More exactly, the rule is the following: to find the first
   appearance of pair \textbf{10}, change it to \textbf{01}
   and put all 1-s to the beginning of string, and all zeros a
   next.  The number of zeros and ones does not change under $P$.  So the orbit of $P$ belongs to the orbit of action of
   infinite symmetric group. Moreover it is clear that the orbit partition of $P$ coincides with the orbit partition of the
   action of the group $S_{\textbf N}$.                                                       \footnote{The initial definition of $P$ starts from the fact
   that the space of paths of infinite Pascal triangle is
   exactly the set of all infinite sequences of 0-1 on
   $\textbf{Z}_2$. In the paper \cite{Ve11} we define  Pascal
   automorphism as the inverse  automorphism $P^{-1}$
   as it was defined above.}

   We have the formulas: $$Px=x+n(x)=T^n(x)x,$$ Remark  that $Px>x$                                                  for all $x$ in the sense of order in integers. So $Px-x$ is a natural (nonnegative) number $n(x)$. The arithmetic
    properties of the (ceiling) function $n: \textbf{Z}_2   \rightarrow \textbf N$ are very interesting.

   The implicit formula for the function $n(\cdot)$ is the
   following:

 $$n(0^m1^k\textbf{10}**)=1^k0^m\textbf{01}-0^m1^k\textbf{10}=
     2^m+2^k-1$$
 Easy calculation shows:
 $$\int_{\textbf{Z}_2}n(x)d\mu=+\infty.$$
 Consequently, the Pascal automorphism {\it is not the  induced
 automorphism by odometer $T$ and vice versa odometer is not
 the integral automorphism over Pascal automorphism.}
 By this construction Pascal automorphism is induced by the
 automorphism $\bar T$ with infinite invariant measure (see next section with the general definitions). Pascal automorphism is a partial case of the general definition of {\it adic}
transformation, see \cite{V82}, \cite{Ve11}].

 The Pascal automorphism as a  transformation on  ${\textbf{Z}_2}$ has the  continuum mutually nonequivalent
 invariant ergodic measures (Bernoulli $(p,q); p+q=1; p,q>0$ measures; this is so called exit boundary of the space of paths of Pascal graph). In the same time the odometer has only one
 ergodic invariant measure (Haar measure). It is possible to
 extent all Bernoulli measures up to invariant measures
 under the automorphism $\bar T$ (see \cite{HIK}).
 There is the sufficiently large literature on Pascal
 automorphism after \cite{V82}(see f.e.\cite{MP,Ve11}, and references therein). The exceptionally interesting papers are \cite{Dl1,Dl2} on the loosely Bernoulli property and
 Takage curve, see recent continuation \cite{LM}.
 There are many open problems about it, e.g. the weak mixing
 (my conjecture in 80-th) is still (11.2015) not proved,
however there are several plans to attack this problem, see \cite{Ve11}. \footnote{The title of that article is not a claim but used to express the old certitude of author that the spectra of Pascal is pure continuous, and the assurance that the proof will be done soon}.

 I will give one example of the problem.
 Important characteristic of Pascal automorphism
  and more general --- adic transformations) is the following.
  Consider  a point
 $x=\{x_n, n \in \textbf N\}\in {\textbf Z}_2$ and for each $n$ fix the $n$-fragment of $x$ of length $n$: $\{x_k, k<n\}$. Suppose that this fragment as a sequence of $0,1$ has $m$ zeros and $n-m$ ones,
 and so belongs to the linear ordered set of all sequences
  of length $n$ and $m$ zeros; suppose $t_n(x)$
 is its ordinal number in this linear ordered set.

   \begin{Hup} (Non-linear limit theorem)
 There exist such a  sequence $a_n$ of natural numbers and  a sequence of  positive numbers $b_n$
 such that for almost all $x\in {\textbf Z}_2$ with respect to
 Haar measure there exists the limit
 $$\lim_n m\{x: \frac{t_n(x)-a_n}{b_n}\leq \alpha \} =\Psi(\alpha),$$
  where $m$ is Haar measure on ${\textbf Z}_2$ and $\Psi(.)$
  is a non-degenerated distribution on $\textbf R$.
  \end{Hup}
   It is intriguing question if this problem has any connection with Takagi curve in the sense of the  paper \cite{Dl2}.

  \section{Generale model}

   \subsection{The notion of generalized induced automorphisms}

  We want to generalize the notion
  of induced automorphism due to S.Kakutani.

  \begin{Th}
  Let $T$ be a m.p. automorphism of the space ($X,\mu)$.
  The automorphism $P$ of the same space is called the generalized induced automorphism over $T$ if one of the following equivalent conditions takes place:

  1. $Px=T^{n(x)}x$, where $n(x)\geq 0$ for almost all $x$.

  2. The orbit partition $\tau (P)$ of automorphism $P$ is a subpartition of the orbit partition of the automorphism $T$:
   $$\tau(P) \succ \tau(T)$$ and the order on the orbits of $P$
   is induced by the order on the orbits of $T$.
   \end{Th}
 The equivalence of the conditions is evident.

 The notion of induced automorphism $T_A$ where $A$ is a measurable set of positive measure  evidently agrees with this generalization. Usually define  $T_Ax=x$ if $x \notin  A$.
 In this case $T_Ax=T^{n(x)}x$ where $n(x)$ is the moment of the first
 return to $A$ if $x \in A$ and $0$ if $x \notin A$; it is clear that
 in this case $\int_X n(x)d\mu <\infty$.

 The class of generalized induced automorphisms is a proper subclass of all automorphisms which has the  form $Px=T^{n(x)}x$
 with arbitrary function $n(\cdot)$ (which defines an automorphism).

 \textbf{Problem.} {\it Describe all automorphisms which are generalized
 induced over the odometer (dyadic or $p$-adic)}.

   It is easy to answer this question in terms of adic transformation:

  \begin{Prop} Suppose the directed graph $\Gamma$ has the property:
  the  number of edges which started from the vertices of the level $n$ to the level $n+1$ does not exceed $p$.
  Then any adic automorphism on the graph is generalized induced over the  $p$-odometer.
   \end{Prop}

   But this description does not answer on the question about the properties of this class of automorphisms like spectra, rank and so on,

   Evidently any automorphism can be represented as generalized induced over an odometers of the spaces $(\prod_n Z_{k_n}, \prod_n m_{k_n})$ with sufficiently large growth of the sequence $k_n$ --- this is a corollary of the universality
   of adic realization (\cite{V82}).

   In the previous examples it is easy to express
   $\int_X n(x)d\mu(x)$ via "binomial coefficients"--- number of paths. Usually it is equal to infinity, so the typical generalized induced automorphism is not ordinary induced.

   \subsection{HK-automorphism as quasi-odometer}

   Pascal automorphism is the simplest generalized induced
   automorphism over dyadic odometer. But because of infinity of the
   integral $\int_X n(x)dm(x)$ it {\it is not the induced
  automorphism by odometer $T$ and consequently odometer
  is not the integral automorphism over Pascal automorphism.}

   In the same time we can define  the integral automorphism over Pascal
   automorphism, using the same function $n(x)$ --- that will be the measure preserving
   automorphism of the space with \textbf{infinite measure}, we call it now  Haijan-Kakutani transformation.

  We have the scheme:

      $$ T\rightarrow P=T^{n(\cdot)} \rightarrow {\bar T}=P^{n(\cdot)}$$

  Here $T$ is 2-odometer, $P$ is Pascal automorphism and $\bar T$ is $HK$ or quasi-odometer which is integral over $P$ and is measure preserving automorphism
  with $\sigma$-finite  measure $\bar \mu$.

  More exactly: define the space with $\sigma$-finite measure and
   measure preserving transformation on it:

     $$(\bar X,\bar\mu, \bar T).$$

   Here $\bar X$ is the  intersection of the lattice $Z^2$ with the subgraph of function $n\geq 0$ (from the formula $Px=T^{n(x)}x$);
the   measure $\bar\mu$ is a $\sigma$ finite measure which coincides with $\mu$ on the base $X$, and whose                                                             conditional measure on the finite sets ($(x,0),(x,1), \dots
  (x,n(x))$ is the counting measure; finally $\bar T$ is defined as  the                                       automorphism which has the automorphism $T$ as induced on the subset $(X,0)$ and acts as $ (x,i)\mapsto (x,i+1)$ if $i=0, \dots n(x)-1$.

    This construction is universal: one can use any automorphism instead of odometer $T$ and any generalized induced over it.
   This gives the class of realizatiosn of automorphisms
   with $\sigma$-finite invariant measure.

   Suppose $T$ and $P$ are m.p.(=measure preserving)
   automorphisms of the space $(X,\mu)$  and $\tau (T),\tau(P)$ are       their orbit partitions; assume that $ \tau(P) \succ \tau (T)$
   which means that each orbit of $P$ belongs to an orbit of
   $T$. In this case  the formula takes place: $$Px=T^{n(x)}x,$$
   where $n(x)$ is a measurable integer-valued function on $x$.
  Assume also that for almost all $x$ the $T$-orbit of
  $x$ consists of infinitely many orbits of $P$}. As was proved in \cite{Ra} this means that the function $n(\cdot)$ is not integrable.
      Remark that this means that the function $n(\cdot)$ is far from
  to be arbitrary integer values function on $X$ --- in contrast to the case of construction of the integral automorphisms where the {\it ceiling function} is arbitrary.
  Thus we have a special class of the integer-valued
  functions $n(\cdot)$ for which we want to give a model
  of infinite measure preserving transformations.

   \begin{Th} Each ergodic automorphism $P$ with finite
   invariant measure can be realized in the form
   $$Px=T^{n(x)}$$ where $T$ is 2-odometer and the orbit partitions $\tau(\cdot)$ of $T$ and $P$ have property
    $ \tau(P) \succ \tau (T)$
    and almost each orbit of $T$ contains countably many
    orbits of $P$.
    \end{Th}
The proof is similar to the proof of Dye's theorem and is based
     on the fact that each ergodic $P$ has
     {\it universal} adic realization in the space of paths
     of the distinguish graph of unordered pairs ($UP$) --- see \cite{VJa} --- this is the strengthen of the theorem on adic
     realization \cite{V82}. But the space of paths of $UP$
     can be identified with the dyadic group $\textbf{Z}_2$ preserving the tail filtrations.

     Of course the dyadic odometer can be changed on any space of
     type  $X=\prod_{k=1}^{\infty}\textbf{p}_k$ where $p_k>1$ is any sequence of naturals, $p_k>1$. Sometimes this is more convenient than $\textbf{p}_n\equiv 2$.

    It is convenient to make the following form for $\bar T$.
    Consider the space $\textbf{Q}_2$ of rational dyadic numbers (as a space but not as an additive group):

    $$ \textbf{Q}_2= \left(\sum_{n<0} {\textbf Z}_2\}\right)\times
    \prod_{n\geq 0}{\textbf Z}_2\equiv \textbf{N}\times               \textbf{Z}_2;$$
    it is useful to consider the first summand as natural
    numbers with addition as operation and second summand as
    dyadic integers with natural 2-adic operation.
    We equip $\textbf{Q}_2$ with (infinite) Haar measure $m$ with normalization  $m(\textbf{Z}_2)=1$.

    Consider the automorphism $P$ which satisfies our assumptions $$Px=T^{n(x)}x=x+n(x), x \in \textbf{Z}_2,$$                     where $T$ is 2-odometer and $n$ is a suitable function.
      {\it Define the automorphism $\bar T$ of  $\textbf{Q}_2$ with
      invariant Haar measure by the following formula:}
        $${\bar T}(y,x)=(y',x'); y,y' \in \sum_{n<0} \textbf{Z}_2,
        x,x' \in  \prod_{n\geq 0} \textbf{Z}_2,$$
        where

        $x'=Px=T^{n(x)}x, \quad \mbox{if} \quad  y=n(x)$;

         $x' =x,  \quad  \mbox{if} \quad  y \ne n(x)$;

         $y'=y+1, \quad \mbox{if} \quad  y \ne n(x)$;

          $y'=0, \quad \mbox{if} \quad  y=n(x)$.

We emphasize that the function $n(x)$ plays two roles -
    as ceiling function and as time change for $P$.
    It is natural to call $\bar T$ the  "pseudo-odometer",
    or "reflection of odometer" with respect to the automorphism $P$.

  {\it It is possible to give another transparent description
  of this construction for $\bar T$:}

 The measure space   can be realized in the space  $\textbf{Q}_2$ in the form which used the previous notations:
  $x \in {\textbf Z}_2; y \in \textbf N, y \leq n(x):$

      $${\bar T}(x,y)=(T^1_y(x),T^2_x(y)),$$

    where $$T^1_x y=y+1, \quad \mbox{if} \quad y<n(x), \quad \mbox{and}\quad T_x^1  y =0 \in \textbf N  \quad \mbox{if}\quad  y=n(x).$$

    For $T^2$ we have:

     $$T^2_y x= x, \quad \mbox{if}\quad y < n(x), T^2_{n(x)}x=Px.$$

     This is simply "horizontal" expression of the fact  $\bar T$ is generalized the induced automorphism. But this form is very convenient for operation and for graph interpretation.

      In general $\{T_x\}$ is a measurable family of permutaions of $\textbf N$ and $\{T_y\}$ is a countable family of ${\textbf Z}_2$   (This construction could be called "bi-skew-product"). But our realization is a little bit special.

    \textbf{Problem} For what ergodic automorphism $\bar T$ with infinite invariant measure there exists such representation (integral representation) where $T$ is odometer and
    with the some function $n(\cdot)$?

    Recall also that the simple generalization of well-known Dye's
    theorem (see f.e.\cite{Krie}) claimed that the orbit partition
    of any ergodic automorphism $S$ is isomorphic to the partition of the space $\textbf{Q}_2$ into  cosets with respect to the subgroup $\textbf Q$  (rational numbers). The model above gives more or less simple way to represent automorphism (like "interval exchange").
      It is possible to give a graph interpretation of
     a representation our model  which is similar to the adic model of the automorphism with
     finite invariant measure.

     For this we need to consider $\textbf Z$-graded (not $\textbf N$=graded) graphs, but of the special type on "minus infinity": the space of double infinite paths must have stabilization on minus infinity,
     which means that there is a path $t$ such that each path coincides with it at $-\infty$.
      We hope that this $\sigma$-finite adic model will be useful for infinite ergodic theory.

\textbf{Problem}
     To define the analog of Bratteli-Vershik model for
     $\sigma$-finite measure preserving transformation.

    Actuality of this problem lies in necessity of the simple models of ergodic automorphisms with $\sigma$-finite measure.

     \section{$\sigma$-finite theory of the random walk on the groups.}

 It was mentioned in the paper \cite{VKa} (see p. 458, item 5.)
 the following important point: the natural theory of the
 random walks on the groups must be developed in the framework of the ergodic theory with infinite invariant measure.  The reason is the possibility to consider a two-sided process and the shift of the trajectories as measure preserving transformation.
 Unfortunately, we are still far from the serious achievements in this direction. The analysis of the new situation must take into account the fact that the global measure is infinite and the space of trajectory of a random walk is not compact. The last facts demand a serious revision of the many ordinary notions
 like ergodicity, regularity, notion of boundary, entropy,
 new look on the role of functional spaces ($l^2,L^{\infty}$), etc.

 Of course the situation seems to be very clear: $\sigma$-finite
 measure in the space of trajectories of random walk is a direct product of Bernoulli measure with  $\sigma$-finite (Haar) measure on the group, and the automorphism is simply skew-product with Bernoulli automorphism as the base and translations on the group as the fibers.
 But this does not mean that all analysis is reduced to the ordinary ergodic theory.

   We give here only evident example of the preference of the "$\sigma$-finite" point of view, more exactly --- preference of the
    consideration of two-sided processes (the time is $\textbf Z$ but not ${\textbf Z}_+$) with necessarily $\sigma$ finite measures.

    How to define exit and entrance (Dynkin)  as well as Poisson-Furstenberg boundaries of the random walk using the shift in the space of trajectories with
 $\sigma$-finite measure?

  Let $G$ be a countable (infinite) finitely generated group and $S$ be the set of generators. We assume that the group $G$ is the semigroup generated by $S$). Denote as $\bar \nu$ the measure on $S^{-1}$:  $(\bar \nu)(g)=\nu(g^{-1})$.

 Consider the random walk on the countable group $G$ with
 a probability measure  $\mu$ on $S$.
  The set of trajectories of generalized Markov chain
  is the subspace $\cal M$ of $G^{\textbf Z}$:
 $${\cal M}=\{\{y_n\}_{n\in {\textbf Z}}: {y_n}^{-1}\cdot y_{n+1}\in S, n \in {\textbf Z}\}.$$

 We can represent the space $\cal M$ as the direct product
 $${\cal M}=S^{\textbf Z}\times G: \{y_n\}=(\{s_n\}_{n\in {\textbf Z}},y_0),$$ where $s_n= {y_n}^{-1}\cdot y_{n+1}, n \in \textbf Z$.
 We will not define a Markov measure on  $\cal M$, but define only transition and co-transition probabilities:
 $$P\{y_n\equiv y_{n-1}s|y_{n-1}\}=\nu(s)\quad
 P\{y_n\equiv y_{n+1}s|y_{n+1}\}=\bar \nu(s)=\nu(s^{-1}),
 n \in \textbf Z.$$

   The dynamics on the space $\cal M$ is defined by the shift $T$
   which is linear {\it automorphism of the space $\cal M$}:
    $\{T(\{y_n\})\}_n=y_{n-1}$.

   We the stationary (shift-invariant) Markov measures on $\cal M$\footnote{Definition of Markov property of the $\sigma$-finite measures is the same as usual:independence of past and future for fixed value of the present.}.

  \begin{Def}
  Exit boundary of a random walk corresponding to the pair $(G,\nu)$  where $G$ is countable group, and  $\nu$ is measure with the finite support $S$ which generate $G$ as a semigroup
  is the set of all ergodic shift-invariant two-sided $\sigma$-finite Markov measures $\mu$ on the space  $\cal M$ with given co-transition probabilities defined above;

  Entrance boundary of random walk corresponding to the pair
  $(G,\mu)$
  is the set of all ergodic shift-invariant two sided $\sigma$-finite
  Markov measures $\mu$ on the space $\cal M$ with transition
  probabilities defined above.
   \end{Def}

  Here we call the Markov measure on the space of double infinite sequences ergodic if it is regular (terms of Kolmogorov) or Kolmogorovian
  (modern term). This means that filtration of the past and filtration of the future have trivial intersections. In other words -no nontrivial events on plus
  and minus infinity.  Remark that the initial distribution does not mentioned in the definition. It is clear that
   both boundaries are topological spaces (with weak topology on the space of measures). It is very important to understand
   the topological properties of these spaces.

  Of course, the exit (correspondingly, entrance) boundary in this
  definition is the same as exit (entrance) boundary  for one-sided system of co-transition (correspondingly, transition) probabilities in the sense of old papers of Dynkin or in the new paper \cite{VJa}. But two-sided invariance gives more possibilities to study these objects.

In the case of symmetric measure $\nu=\bar \nu$ both boundaries are tautologically coincided.

Measure $\nu$ (correspondingly $\bar\nu$) defines a Laplace operators in $l^2(G)$ (these two Laplace operators are mutually conjugate) in the space $l^2$.

    The Laplace operators $L_{\nu}, L_{\bar \nu}$ are convex combinations of the isometries of the shifts $U_s, s\in S$, or $s \in S^{-1}$ in the space $l^2(G)$ (or in the spaces
  $L^{\infty}_{\mu}(\cal M)$).

 The PF-boundary of random walk in this case is a measure space with harmonic measure which is the "part" of the exit
 boundary which is a  topological space. Remark that the exit boundary could be nontrivial even if the PF-boundary is trivial. Both boundaries could be identifying with the spectrum of Laplace
 operator in the space $L^{\infty}$, and more exactly --- with the
 space of minimal non-negative harmonic functions with
 so called harmonic measure.

 One of the most intriguing and (as I know) open questions on
 $\sigma$-finite theory of random walks (even for a  simple walk on $\textbf Z$), is how to describe weakly wandering sets for them.

\kern -0.2cm

\begin{enumerate}
    \item[]{Anatoly Vershik}\\
    {\it St.~Petersburg Department of Stekloff Math. Institute, 			Russian Academy of Sciences\\
    27, Fontanka, 191023\,  St.Petersburg, Russia\\[3pt]}
\end{enumerate}

\end{document}